\documentclass[11pt]{article}
\usepackage{natbib}
\usepackage{amsmath,amsthm,amsopn,amstext,amscd,amsfonts,amssymb,color,enumerate}

\def\tr{\mathop{\rm tr}\nolimits}
\def\erf{\mathop{\rm erf}\nolimits}
\def\Cov{\mathop{\rm Cov}\nolimits}
\def\Var{\mathop{\rm Var}\nolimits}
\def\Vec{\mathop{\rm vec}\nolimits}
\def\Vech{\mathop{\rm vech}\nolimits}
\def\v{\mathop{\rm v}\nolimits}
\def \E {\mathop{\rm E}\nolimits}
\def \P {\mathop{\rm P}\nolimits}
\def \build#1#2#3{\mathrel{\mathop{#1}\limits^{#2}_{#3}}}

\newtheorem{thm}{Theorem}[section]

\newtheorem{lem}{Lemma}[section]
\theoremstyle{definition}

\newtheorem{rem}{Remark}[section]

\newcommand{\ams}[2]
           {\begin{center}
            \begin{minipage}{5.25in}
            \small
            \noindent \textbf{Mathematics Subject Classification: }{\uppercase{#1}}
            \end{minipage}
            \end{center}
            \par\normalsize
           }

\newcommand{\keywords}[1]
           {\begin{center}
            \begin{minipage}{5.25in}
            \small
            \noindent \textbf{Key Words:}~{\textrm{#1}}
            \end{minipage}
            \end{center}
            \normalsize
           }

\title{\vspace{-2.5cm}\textbf{\Large Optimum allocation in multivariate stratified random sampling:
Stochastic matrix optimisation}}
\author{
  \begin{normalsize}
  \begin{tabular}{c}
    \textbf{ Jos\'e A. D\'\i az-Garc\'\i a} \\
    Department of Statistics and Computation \\
    Universidad Aut\'onoma Agraria Antonio Narro \\
    25350 Buenavista, Saltillo, Coahuila, M\'EXICO. \\
    jadiaz@uaaan.mx \\
    and\\
    \textbf{Rogelio Ramos-Quiroga}\\
    Centro de Investigaci\'on en Matem\'aticas\\
    Department of Probability and Statistics\\
    Callej\'on de Jalisco s/n.\\
    36240 Guanajuato, M\'exico\\
    rramosq@cimat.mx\\
  \end{tabular}
  \end{normalsize}
}
\date{}

\begin{document}
  \maketitle
\begin{abstract}
    \noindent The allocation problem for multivariate stratified random sampling as a problem of
    stochastic matrix integer mathematical programming is considered. With these aims the
    asymptotic normality of sample covariance matrices for each strata is established. Some
    alternative approaches are suggested for its solution. An example is solved by applying the
    proposed techniques.
\end{abstract}


\ams{62D05, 90C15, 90C29, 90C10}

\section{\normalsize INTRODUCTION}\label{sec1}

Not long ago, multivariate analysis was mainly based on linear methods illustrated on
small to medium-sized data sets. However, many novel developments, have permitted the
introduction of several innovative statistical and mathematical tools for high-dimensional data
analysis. Developments as generalised multivariate analysis, latent variable analysis, DNA
microarray data, pattern recognition, multivariate  nonlinear analysis, data mining, manifold
learning, shape theory etc., have given a new and modern image to Multivariate Analysis.

One of the topics of statistical theory that is most commonly used in many fields of scientific
research is the theory of probabilistic sampling. From a multivariate point of view, diverse
authors have studied the problem of optimum allocation in multivariate stratified random
sampling. \cite{ad81} and \cite{sssa84}, among many others, proposed the problem of optimum
allocation in multivariate stratified random sampling as a deterministic multiobjective
mathematical programming problem, by considering as objective function a cost function subject
to restrictions on certain functions of variances or viceversa, i.e., considering the functions
of variances as objective and subject to restrictions on costs. Noting that, for the case when
the function of costs is taken as the objective function, the problem of optimum allocation in
multivariate stratified random sampling is reduced to a classical uniobjective mathematical
programming problem. Furthermore, \citet{dgu:08} propose the optimum allocation in multivariate
stratified random sampling as a deterministic nonlinear problem of matrix integer mathematical
programming constrained by a cost function or by a given sample size. Also, \cite{pre78} and
\cite{dggt:07} observe that the values of the population variances are in fact random variables
and formulate the corresponding problem of optimum allocation in multivariate stratified random
sampling as a stochastic mathematical programming problem.

In this paper, the optimum allocation in multivariate stratified random sampling is posed as a
stochastic matrix integer mathematical programming problem constrained by a cost function or by
a given sample size. Section \ref{sec2} provides notation and definitions on multivariate
stratified random sampling. Section \ref{sec3} studies in detail the asymptotic normality of
the sample mean vectors and covariance matrices. The optimum allocation in multivariate
stratified random sampling via stochastic matrix integer mathematical programming is given in
Section \ref{sec4}. Also, several particular solutions are derived for solving the proposed
stochastic mathematical programming problems. Finally, an example of the literature is given in
Section \ref{sec5}.

\section{\normalsize PRELIMINARY RESULTS ON MULTIVARIATE STRATIFIED RANDOM SAMPLING}\label{sec2}

Consider a population of size $N$, divided into $H$ sub-populations (strata). We wish to find a
representative sample of size $n$ and an optimum allocation in the strata meeting the following
requirements: i) to minimise the variance of the estimated mean subject to a budgetary
constraint; or ii) to minimise the cost subject to a constraint on the variances; this is the
classical problem in optimum allocation in univariate stratified sampling, see \cite{coc77},
\cite{sssa84} and \cite{t97}. However, if more than one characteristic (variable) is being considered
then the problem is known as optimum allocation in multivariate stratified sampling. For a
formal expression of the problem of optimum allocation in stratified sampling, consider the
following notation.

The subindex $h=1,2,\cdots,H$ denotes the stratum, $i=1,2,\cdots,N_{h} \mbox{ or } n_{h}$ the unit
within stratum $h$ and $j=1,2,\cdots,G$ denotes the characteristic (variable). Moreover:

\bigskip

\begin{footnotesize}
\begin{tabular}{ll}
    $N_{h}$ & Total number of units within stratum $h$.\\
    $n_{h}$ &  Number of units from the sample in stratum $h$.\\
    \hspace{-.4cm}
    \begin{tabular}{lcl}
       $\mathbf{Y}_{h}$ &=& $(\mathbf{Y}_{h}^{1}, \cdots ,\mathbf{Y}_{h}^{G})$ \\
        &=& $(\mathbf{Y}_{h1}, \cdots ,\mathbf{Y}_{h N_{h}})'$
     \end{tabular} &
     \hspace{-.3cm}
     \begin{tabular}{l}
       $N_{h} \times G$ population matrix in stratum $h$; $\mathbf{Y}_{hi}$ is the\\
       $G$-dimensional value of the $i$-th unit in stratum $h$.\\
     \end{tabular}\\
    \begin{tabular}{lcl}
       $\mathbf{y}_{h}$ &=& $(\mathbf{y}_{h}^{1}, \cdots ,\mathbf{y}_{h}^{G})$ \\
        &=& $(\mathbf{y}_{h1}, \cdots ,\mathbf{y}_{h n_{h}})'$
     \end{tabular} &
     \hspace{-.3cm}
     \begin{tabular}{l}
       $n_{h} \times G$ sample matrix in stratum $h$; $\mathbf{y}_{hi}$ is the $G$-dimensional\\
       $G$-dimensional value of the $i$-th unit of the sample in stratum $h$.\\
     \end{tabular}\\
    $y_{hi}^{j}$ &  Value obtained for the $i$-th unit in stratum $h$\\
    & of the $j$-th characteristic\\[1ex]
    $\mathbf{n} = ({n}_{1},\cdots, {n}_{H})'$ & Vector of the number of units in the sample\\
    $\displaystyle{W_{h}} = \displaystyle{\frac{N_{h}}{N}}$ & Relative size of stratum  $h$\\
    $\displaystyle{\overline{Y}_{h}^{j}} = \frac{1}{N_{h}}\displaystyle
        \sum_{i=1}^{N_{h}}y_{hi}^{j}$ & Population mean in stratum $h$ of the $j$-th characteristic.\\[3ex]
    $\overline{\mathbf{Y}}_{h}=  (\overline{Y}_{h}^{1}, \cdots,\overline{Y}_{h}^{G})'$
    & Population mean vector in stratum $h$. \\[1ex]
    $\displaystyle{\overline{y}_{h}^{j}}=\frac{1}{n_{h}}\displaystyle
        \sum_{i=1}^{n_{h}}y_{hi}^{j}$ & Sample mean in stratum $h$ of the $j$-th characteristic.\\
    $\overline{\mathbf{y}}_{h}=  (\overline{y}_{h}^{1}, \cdots,\overline{y}_{h}^{G})'$
    & Sample mean vector in stratum $h$. \\
\end{tabular}
\end{footnotesize}

\begin{footnotesize}
\begin{tabular}{ll}
    $\displaystyle{\overline{y}_{_{ST}}^{j}=  \sum_{h=1}^{H}W_{h}\overline{y}_{h}^{j}}$
    & Estimator of the population mean in multivariate\\
    & stratified sampling for the $j$-th characteristic.\\
    $\overline{\mathbf{y}}_{_{ST}}=(\overline{y}_{_{ST}}^{1},\cdots,\overline{y}_{_{ST}}^{G})'$
    & Estimator of the population mean vector in\\
    & multivariate stratified sampling.\\
    $\mathbf{S}_{h}$ & Covariance matrix in stratum $h$\\
    & $\mathbf{S}_{h}=\displaystyle\frac{1}{N_{h}}
        \sum_{i=1}^{N_{h}}(\mathbf{y}_{hi}-\overline{\mathbf{Y}}_{h})(\mathbf{y}_{hi}-\overline{\mathbf{Y}}_{h})'$\\[2ex]
    &  where $S_{h_{jk}}$ is the covariance in stratum $h$ of the\\
    & $j$-th and $k$-th characteristics; furthermore\\
    & $S_{h_{jk}}=\displaystyle\frac{1}{{N_{h}}}
        \sum_{i=1}^{N_{h}}(y_{hi}^{j}-\overline{y}_{h}^{j})(y_{hi}^{k}-\overline{y}_{h}^{k})$, and \\
    & $S_{h_{jj}}\equiv S_{hj}^{2} = \displaystyle\frac{1}{N_{h}}
        \sum_{i=1}^{N_{h}}(y_{hi}^{j}-\overline{y}_{h}^{j})^2$.\\
    $\mathbf{s}_{h}$ & Estimator of the covariance matrix in stratum\\
    &  $h$;\\
    & $\mathbf{s}_{h}=\displaystyle\frac{1}{{n_{h}-1}}
        \sum_{i=1}^{n_{h}}(\mathbf{y}_{hi}-\overline{\mathbf{y}}_{h})(\mathbf{y}_{hi}-\overline{\mathbf{y}}_{h})'$\\
    &  where $s_{h_{jk}}$ is the sample covariance in stratum $h$ of the\\
    & $j$-th and $k$-th characteristics; furthermore\\
    & $s_{h_{jk}}=\displaystyle\frac{1}{{n_{h}}-1}
        \sum_{i=1}^{n_{h}}(y_{hi}^{j}-\overline{y}_{h}^{j})(y_{hi}^{k}-\overline{y}_{h}^{k})$, and \\
    & $s_{h_{jj}}\equiv s_{hj}^{2} = \displaystyle\frac{1}{{n_{h}}-1}
        \sum_{i=1}^{n_{h}}(y_{hi}^{j}-\overline{y}_{h}^{j})^2.$\\
    $\Cov(\overline{\mathbf{y}}_{_{ST}})$
        & Covariance matrix of $\overline{\mathbf{y}}_{_{ST}}.$\\
    $
    \widehat{\Cov}(\overline{\mathbf{y}}_{_{ST}})$
        & Estimator of the covariance matrix of $\overline{\mathbf{y}}_{_{ST}}$,\\
        & it is denoted as $\widehat{\Cov}(\overline{\mathbf{y}}_{_{ST}}) \equiv
        \widehat{\Cov(\overline{\mathbf{y}}_{_{ST}})}$, and defined as \\[1ex]
     & $
             = \left(
                \begin{array}{cccc}
                  \widehat{\Var}(\overline{y}_{_{ST}}^{1}) & \widehat{\Cov}(\overline{y}_{_{ST}}^{1},\overline{y}_{_{ST}}^{2})
                    & \cdots & \widehat{\Cov}(\overline{y}_{_{ST}}^{1},\overline{y}_{_{ST}}^{G}) \\
                  \widehat{\Cov}(y_{_{ST}}^{2},\overline{y}_{_{ST}}^{1}) & \widehat{\Var}(\overline{y}_{_{ST}}^{2}) & \cdots
                    & \widehat{\Cov}(\overline{y}_{_{ST}}^{2},\overline{y}_{_{ST}}^{G}) \\
                  \vdots & \vdots & \ddots & \vdots \\
                  \widehat{\Cov}(\overline{y}_{_{ST}}^{G},\overline{y}_{_{ST}}^{1}) & \widehat{\Cov}(\overline{y}_{_{ST}}^{G},
                    \overline{y}_{_{ST}}^{2}) & \cdots & \widehat{\Var}(\overline{y}_{_{ST}}^{G}) \\
                \end{array}
            \right )
          $\\[3ex]
        & = $\displaystyle{\sum_{h=1}^{H}\frac{{{W_{h}}^{2}}
        \mathbf{s}_{h}}{n_{h}} - \sum_{h=1}^{H} \frac{{W_{h}}\mathbf{s}_{h}}{N}}$\\[1ex]
        $\widehat{\Cov}(\overline{y}_{_{ST}}^{j},\overline{y}_{_{ST}}^{k}) $ &
        Estimated covariance of $\overline{y}_{_{ST}}^{j}$ and $\overline{y}_{_{ST}}^{k}$ where \\[1ex]
        & $\widehat{\Cov}(\overline{y}_{_{ST}}^{k},\overline{y}_{_{ST}}^{j}) \equiv
        \widehat{\Cov(\overline{y}_{_{ST}}^{j},\overline{y}_{_{ST}}^{k})}$, with \\[1ex]
        \hspace{3.5cm}
        &  $\widehat{\Cov}(\overline{y}_{_{ST}}^{j},\overline{y}_{_{ST}}^{k})= \displaystyle{\sum_{h=1}^{H}\frac{{{W_{h}}^{2}}
        s_{h_{jk}}}{n_{h}} - \sum_{h=1}^{H} \frac{{W_{h}}s_{h_{jk}}}{N}}$, and \\[1ex]
        &
        $\widehat{\Cov}(\overline{y}_{_{ST}}^{j},\overline{y}_{_{ST}}^{j}) \equiv
        \widehat{\Var}(\overline{y}_{_{ST}}^{j})= \displaystyle{\sum_{h=1}^{H}\frac{{{W_{h}}^{2}}
        s_{hj}^{2}}{n_{h}} - \sum_{h=1}^{H} \frac{{W_{h}}s_{hj}^{2}}{N}}$. \\[2ex]
       $c_{h}$ & Cost per $G$-dimensional sampling unit in stratum $h$ and let\\
       & $\mathbf{c} = (c_{1}, \dots, c_{G})'$.\\
\end{tabular}
\end{footnotesize}

\noindent Where if $\mathbf{a} \in \Re^{G}$, $\mathbf{a}'$ denotes the transpose of
$\mathbf{a}$.

\section{\normalsize LIMITING DISTRIBUTION OF SAMPLE MEANS AND COVARIANCE MATRICES}\label{sec3}

In this section the asymptotic distribution of the estimator of the covariance
matrix $\mathbf{s}_{h}$ and mean $\overline{\mathbf{y}}_{h}$  is considered. With this aim in mind,
the multivariate version of H\'ajek's theorem is proposed in the context of sampling theory in
terms of the extension stated in \citet{h:61}. First, consider the following notation and
definitions.

A detailed discussion of operator ``$\Vec$", ``$\Vech$", Moore-Penrose inverse, Kronecker
product, commutation matrix and duplication matrix may be found in \citet{mn:88}, among many
others. For convenience, some notations shall be introduced, although in general it adheres to
standard notations.

For all matrix $\mathbf{A}$, there exists a unique matrix $\mathbf{A}^{+}$ which is termed
the \emph{Moore-Penrose inverse} of $\mathbf{A}$.

Let $\mathbf{A}$ be an $m \times n$ matrix and $\mathbf{B}$ a $p \times q$ matrix. The $mp
\times nq$ matrix defined by
$$
  \left[
  \begin{array}{ccc}
    a_{11}\mathbf{B} & \cdots & a_{11}\mathbf{B} \\
    \vdots & \ddots & \vdots \\
    a_{11}\mathbf{B} & \cdots & a_{11}\mathbf{B}
  \end{array}
  \right]
$$
is termed the \emph{Kronecker product} (also termed tensor product or direct product) of
$\mathbf{A}$ and $\mathbf{B}$ and written $\mathbf{A} \otimes \mathbf{B}$. Let $\mathbf{C}$ be
an $m \times n$ matrix and $\mathbf{C}_{j}$ its $j$-th column, then $\Vec \mathbf{C}$ is the
$mn \times 1$ vector
$$
  \Vec \mathbf{C} =
  \left [
  \begin{array}{c}
    \mathbf{C}_{1} \\
    \mathbf{C}_{2} \\
    \vdots \\
    \mathbf{C}_{n}
  \end{array}
  \right].
$$
The vector $\Vec \mathbf{C}$ and $\Vec \mathbf{C}^{'}$ clearly contain the same $mn$
components, but in different order. Therefore there exist a unique $mn \times mn$ permutation
matrix which transform $\Vec \mathbf{C}$ into $\Vec \mathbf{C}'$. This matrix is termed the
\emph{commutation matrix} and is denoted $\mathbf{K}_{mn}$. (If $m=n$, is often write
$\mathbf{K}_{n}$ instead of $\mathbf{K}_{mn}$.) Hence
$$
  \mathbf{K}_{mn} \Vec \mathbf{C} = \Vec \mathbf{C}'.
$$
Similarly, let $\mathbf{B}$ be a square $n \times n$ matrix. Then $\Vech \mathbf{B}$ (also
denoted as $\v(\mathbf{B})$) shall denote the $n(n+1)/2 \times 1$ vector that is obtained
from $\Vec \mathbf{B}$ by eliminating all supradiagonal elements of $\mathbf{B}$. If
$\mathbf{B} = \mathbf{B}'$, $\Vech \mathbf{B}$ contains only the distinct elements of
$\mathbf{B}$, then there is a unique $n^{2} \times n(n+1)/2$ matrix termed \textit{duplication
matrix}, which is denoted by $\mathbf{D}_{n}$, such that $\mathbf{D}_{n}\Vech \mathbf{B} = \Vec
\mathbf{B}$ and $\mathbf{D}_{n}^{+}\Vec \mathbf{B} = \Vech \mathbf{B}$. Finally, denote
$(\Vech \mathbf{B})' \equiv \Vech' \mathbf{B}$.

In what follows, from Lemma \ref{lemma1} through Theorem \ref{teo2}, asymptotic results are stated
for a single stratum. The notation $N_{\nu}$ and $n_{\nu}$ denote the size of a generic stratum
and the size of a simple random sample from that stratum.

\begin{lem}\label{lemma1}
Let $\mathbf{\mathbf{\Xi}}_{\nu}$ be a $G \times G$ symmetric random matrix defined as
$$
   \mathbf{\Xi}_{\nu} =  \frac{1}{n_{\nu}-1}\sum_{i = 1}^{n_{\nu}}(\mathbf{y}_{\nu i}- \overline{\mathbf{Y}}_{\nu})
   (\mathbf{y}_{\nu i}- \overline{\mathbf{Y}}_{\nu})'.
$$
Suppose that for $\boldsymbol{\lambda} = (\lambda_{1}, \dots, \lambda_{k})'$, any vector of constants, $k = G(G+1)/2$,%
{\small
\begin{equation}\label{shc}
  \boldsymbol{\lambda}'\left(\mathbf{M}_{\nu}^{4} - \Vech \mathbf{S}_{\nu}
  \Vech' \mathbf{S}_{\nu}\right) \boldsymbol{\lambda} \geq \epsilon \build{\max}{}{1 \leq \alpha \leq k
  } \left[\lambda_{\alpha}^{2} \mathbf{e}_{k}^{\alpha '}\left(\mathbf{M}_{\nu}^{4} - \Vech \mathbf{S}_{\nu}
  \Vech' \mathbf{S}_{\nu}\right) \mathbf{e}_{k}^{\alpha }\right],
\end{equation}}
where $\mathbf{e}_{k}^{\alpha } = (0, \dots, 0, 1, 0, \dots, 0)'$ is the $\alpha$-th vector of
the canonical base of $\Re^{k}$, $\epsilon > 0$ and independent of $\nu
> 1$ and
$$
  \mathbf{M}_{\nu}^{4} = \frac{1}{N_{\nu}}\mathbf{D}_{G}^{+}\left[\sum_{i = 1}^{N_{\nu}}
  (\mathbf{y}_{\nu i} - \overline{\mathbf{Y}}_{\nu})(\mathbf{y}_{\nu i} - \overline{\mathbf{Y}}_{\nu})'
  \otimes (\mathbf{y}_{\nu i} - \overline{\mathbf{Y}}_{\nu})(\mathbf{y}_{\nu i} - \overline{\mathbf{Y}}_{\nu})'
  \right]\mathbf{D}_{G}^{+'},
$$
is the fourth central moment. Assume that $n_{\nu}\rightarrow \infty$, $N_{\nu} - n_{\nu}
\rightarrow \infty$, $N_{\nu}\rightarrow \infty$, and that, for all $j = 1,\dots,G$,%
\begin{equation}\label{hcas}
  \left[\build{\lim}{}{\nu \rightarrow \infty}\left(\frac{n_{\nu}}{N_{\nu}}\right) = 0\right] \Rightarrow
  \build{\lim}{}{\nu \rightarrow \infty} \frac{\build{\max}{}{1 \leq i_{1} < \cdots < i_{n_{\nu}}\leq N_{\nu}}
  \displaystyle\sum_{\beta = 1}^{n_{\nu}}\left[\left(y_{\nu i_{\beta}}^{j} - \overline{Y}_{\nu}^{j}\right)^{2} - S_{\nu j}^{2}\right]^{2}}
  {N_{\nu}\left[m_{\nu j}^{4}-\left(S_{\nu j}^{2}\right)^{2}\right]} = 0,
\end{equation}
where
$$
  m_{\nu j}^{4} =\frac{1}{N_{\nu}}\sum_{i = 1}^{N_{\nu}}\left(y_{\nu i}^{j} -
  \overline{y}_{\nu}^{j}\right)^{4}.
$$
Then, $\Vech\mathbf{\Xi}_{\nu}$ is asymptotically normally distributed as
$$
  \Vech\mathbf{\Xi}_{\nu} \build{\rightarrow}{d}{} \mathcal{N}_{k}(\E(\Vech\mathbf{\Xi}_{\nu}),
  \Cov(\Vech\mathbf{\Xi}_{\nu})),
$$
with
\begin{equation}\label{mXi}
    \E(\Vech\mathbf{\Xi}_{\nu}) = \frac{n_{\nu}}{n_{\nu}-1}\Vech\mathbf{S}_{\nu},
\end{equation}
and
\begin{equation}\label{cmXi}
    \Cov(\Vech\mathbf{\Xi}_{\nu}) = \frac{n_{\nu}}{(n_{\nu} - 1)^{2}}\left(\mathbf{M}_{\nu}^{4}
     - \Vech \mathbf{S}_{\nu}\Vech' \mathbf{S}_{\nu}\right).
\end{equation}
$n_{\nu}$ is the sample size for a simple random sample from the $\nu$-th population of size
$N_{\nu}$.
\end{lem}

\begin{rem}\label{remark1}
Let
$$
   \mathbf{\Xi}_{\nu} =  \frac{1}{n_{\nu}-1}\sum_{i = 1}^{n_{\nu}}(\mathbf{y}_{\nu i}- \overline{\mathbf{Y}}_{\nu})
   (\mathbf{y}_{\nu i}- \overline{\mathbf{Y}}_{\nu})'.
$$
Hence,
\begin{eqnarray*}
  \Vec \mathbf{\Xi}_{\nu} &=&  \frac{1}{n_{\nu}-1}\sum_{i = 1}^{n_{\nu}}\Vec(\mathbf{y}_{\nu i}- \overline{\mathbf{Y}}_{\nu})
   (\mathbf{y}_{\nu i}- \overline{\mathbf{Y}}_{\nu})' \\
   &=&  \frac{1}{n_{\nu}-1}\sum_{i = 1}^{n_{\nu}}(\mathbf{y}_{\nu i}- \overline{\mathbf{Y}}_{\nu})
   \otimes (\mathbf{y}_{\nu i}- \overline{\mathbf{Y}}_{\nu}).
\end{eqnarray*}
From where
$$
  \Vech \mathbf{\Xi}_{\nu} =  \frac{1}{n_{\nu}-1}\sum_{i = 1}^{n_{\nu}}\mathbf{D}_{G}^{+}(\mathbf{y}_{\nu i}- \overline{\mathbf{Y}}_{\nu})
   \otimes (\mathbf{y}_{\nu i}- \overline{\mathbf{Y}}_{\nu}),
$$
$k = G(G+1)/2$.

Taking $m = k$ and $\mathbf{a}_{\nu i} = (a_{\nu i}^{1}, \dots, a_{\nu i}^{k})' =
\mathbf{D}_{G}^{+}(\mathbf{y}_{\nu i}- \overline{\mathbf{Y}}_{\nu}) \otimes (\mathbf{y}_{\nu
i}- \overline{\mathbf{Y}}_{\nu})$ in \citet{h:61}, it is obtained that:
\begin{enumerate}[i)]
  \item $\Vech \mathbf{\Xi}_{\nu}$ can be expressed as
  $$
    \Vech \mathbf{\Xi}_{\nu} = \sum_{i =1}^{N_{\nu}} b_{\nu i} \mathbf{a}_{\nu R_{\nu i}}.
  $$
  with $b's$ fixed, furthermore $b_{\nu 1} = \cdots = b_{\nu n_{\nu}} = 1/(n_{\nu}-1)$, $b_{\nu n_{\nu}+1} = \cdots = b_{\nu N_{\nu}} =
  0$. Then
  $$
    \build{\lim}{}{\nu \rightarrow \infty} \frac{\build{\max}{}{1 \leq j \leq N_{\nu}} \left(b_{\nu j} -
    \overline{b}_{\nu}\right)^{2}}{\displaystyle\sum_{i=1}^{N_{\nu}}\left(b_{\nu j} -
    \overline{b}_{\nu}\right)^{2}} = 0, \quad \mbox{ where } \quad \overline{b}_{\nu} = \frac{1}{N_{\nu}}\sum_{i = 1}^{N_{\nu}}b_{\nu i}
  $$
  holds if $n_{\nu}\rightarrow \infty$, $N_{\nu} - n_{\nu} \rightarrow \infty$.

  \item $\overline{\mathbf{a}}_{\nu} = (\overline{a}_{\nu}^{1} \cdots \overline{a}_{\nu}^{k})'$ is
  \begin{eqnarray*}
    \overline{\mathbf{a}}_{\nu} &=& \frac{1}{N_{\nu}} \sum_{i = 1}^{N_{\nu}}\mathbf{a}_{\nu i}\\
     &=& \frac{1}{N_{\nu}} \sum_{i = 1}^{N_{\nu}} \mathbf{D}_{G}^{+} (\mathbf{y}_{\nu i}- \overline{\mathbf{Y}}_{\nu})
     \otimes (\mathbf{y}_{\nu i}- \overline{\mathbf{Y}}_{\nu})\\
     &=& \Vech \frac{1}{N_{\nu}} \sum_{i = 1}^{N_{\nu}} (\mathbf{y}_{\nu i}- \overline{\mathbf{Y}}_{\nu})
     (\mathbf{y}_{\nu i}- \overline{\mathbf{Y}}_{\nu})'\\
     &=& \Vech \mathbf{S}_{\nu}
  \end{eqnarray*}
  \item From (7.2) in \citet{h:61}
  \begin{equation}\label{ohc}
    \sum_{i = 1}^{N_{\nu}}\left[\sum_{\alpha = 1}^{k}\lambda_{\alpha} (a_{\nu i}^{\alpha} -
    a_{\nu}^{\alpha})\right]^{2} \geq \epsilon \build{\max}{}{1 \leq \alpha \leq k}
    \left[\lambda_{\alpha}^{2}\sum_{i = 1}^{N_{\nu}} (a_{\nu i}^{\alpha} -
    a_{\nu}^{\alpha})^{2}\right].
  \end{equation}
  In the context of sampling theory the right side in (\ref{ohc}) can be written as
  \begin{eqnarray*}
  \hspace{-1cm}
    \sum_{i = 1}^{N_{\nu}}\left[\sum_{\alpha = 1}^{k}\lambda_{\alpha} (a_{\nu i}^{\alpha} - a_{\nu}^{\alpha})\right]^{2}
     &=& \sum_{i = 1}^{N_{\nu}}\left\{\boldsymbol{\lambda}' \left[\mathbf{D}_{G}^{+} (\mathbf{y}_{\nu i}- \overline{\mathbf{Y}}_{\nu})
     \otimes (\mathbf{y}_{\nu i}- \overline{\mathbf{Y}}_{\nu}) - \Vech \mathbf{S}_{\nu}\right]\right\}^{2} \\
  \end{eqnarray*}

  \vspace{-1.5cm}

  \begin{eqnarray}
     \quad &=& \sum_{i = 1}^{N_{\nu}}\boldsymbol{\lambda}' \left[\mathbf{D}_{G}^{+} (\mathbf{y}_{\nu i}- \overline{\mathbf{Y}}_{\nu})
     \otimes (\mathbf{y}_{\nu i}- \overline{\mathbf{Y}}_{\nu}) - \Vech \mathbf{S}_{\nu}\right]\nonumber\\
     &&  \left[(\mathbf{y}_{\nu i}- \overline{\mathbf{Y}}_{\nu})' \otimes (\mathbf{y}_{\nu i}- \overline{\mathbf{Y}}_{\nu})'\mathbf{D}_{G}^{+'} -
     \Vech' \mathbf{S}_{\nu}\right]\boldsymbol{\lambda}\nonumber\\
     &=& \boldsymbol{\lambda}' \left[\mathbf{D}_{G}^{+} \sum_{i = 1}^{N_{\nu}}(\mathbf{y}_{\nu i}- \overline{\mathbf{Y}}_{\nu})
     (\mathbf{y}_{\nu i}- \overline{\mathbf{Y}}_{\nu})' \otimes (\mathbf{y}_{\nu i}- \overline{\mathbf{Y}}_{\nu})
     (\mathbf{y}_{\nu i}- \overline{\mathbf{Y}}_{\nu})'\mathbf{D}_{G}^{+'}\right.\nonumber\\
     && - \Vech \mathbf{S}_{\nu} \sum_{i = 1}^{N_{\nu}}(\mathbf{y}_{\nu i}-
     \overline{\mathbf{Y}}_{\nu})' \otimes (\mathbf{y}_{\nu i}- \overline{\mathbf{Y}}_{\nu})'\mathbf{D}_{G}^{+'}\nonumber\\
     && \left. - \mathbf{D}_{G}^{+}\sum_{i = 1}^{N_{\nu}}(\mathbf{y}_{\nu i}-
     \overline{\mathbf{Y}}_{\nu}) \otimes (\mathbf{y}_{\nu i}- \overline{\mathbf{Y}}_{\nu})\Vech' \mathbf{S}_{\nu}
     + N_{\nu}\Vech \mathbf{S}_{\nu}\Vech'
     \mathbf{S}_{\nu}\right]\boldsymbol{\lambda}\nonumber\\ \label{lshc}
     &=& N_{\nu}\boldsymbol{\lambda}'\left(\mathbf{M}_{\nu}^{4} - \Vech \mathbf{S}_{\nu}
     \Vech' \mathbf{S}_{\nu}\right) \boldsymbol{\lambda},
  \end{eqnarray}
  where $\mathbf{M}_{\nu}^{4}$ is%
  {\small
  \begin{equation}\label{m4}
     = \frac{1}{N_{\nu}}\mathbf{D}_{G}^{+}\left[\sum_{i = 1}^{N_{\nu}}
    (\mathbf{y}_{\nu i} - \overline{\mathbf{Y}}_{\nu})(\mathbf{y}_{\nu i} - \overline{\mathbf{Y}}_{\nu})'
    \otimes (\mathbf{y}_{\nu i} - \overline{\mathbf{Y}}_{\nu})(\mathbf{y}_{\nu i} - \overline{\mathbf{Y}}_{\nu})'
    \right]\mathbf{D}_{G}^{+'},
  \end{equation}}
  Similarly the right side of (\ref{ohc}) is
  \begin{eqnarray*}
  \hspace{-1cm}
    \lambda_{\alpha}^{2}\sum_{i = 1}^{N_{\nu}} (a_{\nu i}^{\alpha} - a_{\nu}^{\alpha})^{2}
     &=& \sum_{i = 1}^{N_{\nu}}\left\{\boldsymbol{\lambda}' \mathbf{e}_{k}^{\alpha}\mathbf{e}_{k}^{\alpha'}
     \left[\mathbf{D}_{G}^{+} (\mathbf{y}_{\nu i}- \overline{\mathbf{Y}}_{\nu})
     \otimes (\mathbf{y}_{\nu i}- \overline{\mathbf{Y}}_{\nu}) - \Vech \mathbf{S}_{\nu}\right]\right\}^{2} \\
     &=& \lambda_{\alpha}^{2}\sum_{i = 1}^{N_{\nu}}\left\{\mathbf{e}_{k}^{\alpha'}
     \left[\mathbf{D}_{G}^{+} (\mathbf{y}_{\nu i}- \overline{\mathbf{Y}}_{\nu})
     \otimes (\mathbf{y}_{\nu i}- \overline{\mathbf{Y}}_{\nu}) - \Vech
     \mathbf{S}_{\nu}\right]\right\}^{2}.
  \end{eqnarray*}
  Then, proceeding as in 3.,
  \begin{equation}\label{rshc}
    \lambda_{\alpha}^{2}\sum_{i = 1}^{N_{\nu}} (a_{\nu i}^{\alpha} - a_{\nu}^{\alpha})^{2} =
    N_{\nu} \lambda_{\alpha}^{2} \mathbf{e}_{k}^{\alpha '}\left(\mathbf{M}_{\nu}^{4} -
    \Vech \mathbf{S}_{\nu}\Vech' \mathbf{S}_{\nu}\right) \mathbf{e}_{k}^{\alpha}.
  \end{equation}
  Therefore, from (\ref{lshc}) and (\ref{rshc}), (\ref{shc}) is established.
  \item The expression for (\ref{hcas}) is found analogously as the procedure described in item 3.
  \item Finally,
  \begin{eqnarray*}
    \E(\Vech \mathbf{\Xi}) &=& \frac{1}{n_{\nu}-1} \sum_{i = 1}^{n_{\nu}}\E\mathbf{D}_{G}^{+}(\mathbf{y}_{\nu i}-
    \overline{\mathbf{Y}}_{\nu}) \otimes (\mathbf{y}_{\nu i}- \overline{\mathbf{Y}}_{\nu}) \\
     &=& \frac{1}{n_{\nu}-1} \sum_{i = 1}^{n_{\nu}}\Vech \E(\mathbf{y}_{\nu i}-
    \overline{\mathbf{Y}}_{\nu})(\mathbf{y}_{\nu i}- \overline{\mathbf{Y}}_{\nu})' \\
     &=& \frac{1}{n_{\nu}-1} \sum_{i = 1}^{n_{\nu}} \Vech \mathbf{S}_{\nu}\\
     &=& \frac{n_{\nu}}{n_{\nu}-1}  \Vech \mathbf{S}_{\nu}\\
  \end{eqnarray*}
  Similarly, by independence
  $$\hspace{-2cm}
    \Cov(\Vech \mathbf{\Xi}) = \frac{1}{(n_{\nu}-1)^{2}} \sum_{i = 1}^{n_{\nu}}\Cov\left[\mathbf{D}_{G}^{+}(\mathbf{y}_{\nu i}-
    \overline{\mathbf{Y}}_{\nu}) \otimes (\mathbf{y}_{\nu i}- \overline{\mathbf{Y}}_{\nu})\right]
  $$

  \vspace{-.75cm}

  {\small
  \begin{eqnarray*}
     \quad &=& \frac{1}{(n_{\nu}-1)^{2}} \sum_{i = 1}^{n_{\nu}}\left\{\E\left[\mathbf{D}_{G}^{+}(\mathbf{y}_{\nu i}-
    \overline{\mathbf{Y}}_{\nu}) \otimes (\mathbf{y}_{\nu i}- \overline{\mathbf{Y}}_{\nu})(\mathbf{y}_{\nu i}-
    \overline{\mathbf{Y}}_{\nu})' \otimes (\mathbf{y}_{\nu i}- \overline{\mathbf{Y}}_{\nu})'\mathbf{D}_{G}^{+'}\right]
    \right .\\
     && \left. - \E\left[\mathbf{D}_{G}^{+}(\mathbf{y}_{\nu i}- \overline{\mathbf{Y}}_{\nu}) \otimes
     (\mathbf{y}_{\nu i}- \overline{\mathbf{Y}}_{\nu})\right] \E \left[(\mathbf{y}_{\nu i}- \overline{\mathbf{Y}}_{\nu})'
     \otimes (\mathbf{y}_{\nu i}- \overline{\mathbf{Y}}_{\nu})'\mathbf{D}_{G}^{+'}\right] \right\}\\
     &=& \frac{1}{(n_{\nu}-1)^{2}} \sum_{i = 1}^{n_{\nu}} \left(\mathbf{M}_{\nu}^{4} - \Vech \mathbf{S}_{\nu}\Vech' \mathbf{S}_{\nu}\right)\\
     &=& \frac{n_{\nu}}{(n_{\nu}-1)^{2}}  \left(\mathbf{M}_{\nu}^{4} - \Vech
     \mathbf{S}_{\nu}\Vech' \mathbf{S}_{\nu}\right),
  \end{eqnarray*}}
  the last expression is obtained observing that
  $$
    \E\left[\mathbf{D}_{G}^{+}(\mathbf{y}_{\nu i}- \overline{\mathbf{Y}}_{\nu}) \otimes
     (\mathbf{y}_{\nu i}- \overline{\mathbf{Y}}_{\nu})\right] = \Vech \E\left[(\mathbf{y}_{\nu i}- \overline{\mathbf{Y}}_{\nu})
     (\mathbf{y}_{\nu i}- \overline{\mathbf{Y}}_{\nu})'\right] = \Vech \mathbf{S}_{\nu}
  $$
  and that
  $$
    \E\left\{\mathbf{D}_{G}^{+}
    (\mathbf{y}_{\nu i} - \overline{\mathbf{Y}}_{\nu})(\mathbf{y}_{\nu i} - \overline{\mathbf{Y}}_{\nu})'
    \otimes (\mathbf{y}_{\nu i} - \overline{\mathbf{Y}}_{\nu})(\mathbf{y}_{\nu i} - \overline{\mathbf{Y}}_{\nu})'
    \mathbf{D}_{G}^{+'}\right\} = \mathbf{M}_{\nu}^{4}
  $$
  where $\mathbf{M}_{\nu}^{4}$ is defined in (\ref{m4}). \qed
\end{enumerate}
\end{rem}

\begin{thm}\label{teo1}
Under assumptions in Lemma \ref{lemma1}, the sequence of sample covariance matrices
$\mathbf{s}_{\nu}$ are such that $\Vech \mathbf{s}_{\nu}$ has an asymptotic normal distribution with
asymptotic mean and covariance matrix given by (\ref{mXi}) and (\ref{cmXi}), respectively.
\end{thm}
\begin{proof}
This follows immediately from Lemma \ref{lemma1}, only observe that
\begin{eqnarray*}
  \mathbf{s}_{\nu} &=&  \frac{1}{n_{\nu}-1}\sum_{i = 1}^{n_{\nu}}(\mathbf{y}_{\nu i}- \overline{\mathbf{y}}_{\nu})
   (\mathbf{y}_{\nu i}- \overline{\mathbf{y}}_{\nu})' \\
   &=& \mathbf{\Xi} - \frac{n_{\nu}}{n_{\nu}-1} (\overline{\mathbf{y}}_{\nu}- \overline{\mathbf{Y}}_{\nu})
   (\overline{\mathbf{y}}_{\nu }- \overline{\mathbf{Y}}_{\nu})',
\end{eqnarray*}
where
$$
  \frac{n_{\nu}}{n_{\nu}-1} \rightarrow 1 \quad \mbox{ and } \quad (\overline{\mathbf{y}}_{\nu}- \overline{\mathbf{Y}}_{\nu})
   (\overline{\mathbf{y}}_{\nu }- \overline{\mathbf{Y}}_{\nu})'\rightarrow 0 \quad \mbox{in
   probability}. \qquad\mbox{\qed}
$$
\end{proof}

\begin{rem}
Observe that it is possible to find the asymptotic distribution of $\Vec \mathbf{s}_{\nu}$, but
this asymptotic normal distribution is singular, because $\Cov(\Vec\mathbf{s}_{\nu})$ is
singular. This is due to the fact $\Cov(\Vec\mathbf{s}_{\nu})$ is the $G^{2} \times G^{2}$
covariance matrix in the asymptotic distribution distribution of $\Vec \mathbf{s}_{\nu}$ and,
because $\mathbf{s}_{\nu}$ is symmetric, then $\Vec \mathbf{s}_{\nu}$ has repeated elements. In
this case, $\Vec\mathbf{s}_{\nu}$ is asymptotically normally distributed as (see \citet{mh:82})
$$
  \Vec\mathbf{s}_{\nu} \build{\rightarrow}{d}{} \mathcal{N}_{G^{2}}(\E(\Vec\mathbf{\Xi}_{\nu}),
  \Cov(\Vec\mathbf{\Xi}_{\nu})),
$$
where
$$
  \E(\Vec\mathbf{\Xi}_{\nu}) = \frac{n_{\nu}}{n_{\nu}-1}\Vec\mathbf{S}_{\nu},
$$
$$
  \Cov(\Vec\mathbf{\Xi}_{\nu}) = \frac{n_{\nu}}{(n_{\nu} - 1)^{2}}\left(\mathfrak{M}_{\nu}^{4}
  - \Vec \mathbf{S}_{\nu}\Vec' \mathbf{S}_{\nu}\right),
$$
and
$$
  \mathfrak{M}_{\nu}^{4} = \frac{1}{N_{\nu}}\left[\sum_{i = 1}^{N_{\nu}}
  (\mathbf{y}_{\nu i} - \overline{\mathbf{Y}}_{\nu})(\mathbf{y}_{\nu i} - \overline{\mathbf{Y}}_{\nu})'
  \otimes (\mathbf{y}_{\nu i} - \overline{\mathbf{Y}}_{\nu})(\mathbf{y}_{\nu i} - \overline{\mathbf{Y}}_{\nu})'
  \right]. \quad\mbox{\qed}
$$
\end{rem}

Proceeding in analogous way as in Lemma \ref{lemma1} and  Remark \ref{remark1}, it is
obtained:

\begin{thm}\label{teo2}
Suppose that for $\boldsymbol{\lambda} = (\lambda_{1}, \dots, \lambda_{G})'$, any vector of
constants,
\begin{equation}\label{shcm}
  \boldsymbol{\lambda}'\mathbf{S}_{\nu}\boldsymbol{\lambda} \geq \epsilon \build{\max}{}{1 \leq j
  \leq G} \left[\lambda_{\alpha}^{2} S_{\nu \alpha}^{2}\right].
\end{equation}
Assume that $n_{\nu}\rightarrow \infty$, $N_{\nu} - n_{\nu}
\rightarrow \infty$, $N_{\nu}\rightarrow \infty$, and that%
\begin{equation}\label{hcas1}
  \left[\build{\lim}{}{\nu \rightarrow \infty}\left(\frac{n_{\nu}}{N_{\nu}}\right) = 0\right] \Rightarrow
  \build{\lim}{}{\nu \rightarrow \infty} \frac{\build{\max}{}{1 \leq i_{1} < \cdots < i_{n_{\nu}}\leq N_{\nu}}
  \displaystyle\sum_{\beta = 1}^{n_{\nu}}\left(y_{\nu i_{\beta}}^{j} - \overline{Y}_{\nu}^{j}\right)^{2}}
  {N_{\nu}S_{\nu j}^{2}} = 0,
\end{equation}
Then, $\overline{\mathbf{y}}_{\nu}$ is asymptotically normally distributed as
$$
  \overline{\mathbf{y}}_{\nu} \build{\rightarrow}{d}{} \mathcal{N}_{G}\left(\overline{\mathbf{Y}}_{\nu},
  \mathbf{S}_{\nu}\right).
$$
$n_{\nu}$ is the sample size for a simple random sample from the $\nu$-th population of size
$N_{\nu}$.
\end{thm}

As direct consequence of Theorem \ref{teo1} it is obtained:

\begin{thm}
Let $\widehat{\Cov}(\overline{\mathbf{y}}_{_{ST}})$ be the estimator of the covariance matrix
of $\overline{\mathbf{y}}_{ST}$, then
$$
  \Vech \widehat{\Cov}(\overline{\mathbf{y}}_{_{ST}}) = \sum_{h=1}^{H}\left(\frac{{{W_{h}}^{2}}}{n_{h}} -
  \frac{{W_{h}}}{N} \right)\Vech \mathbf{s}_{h}
$$
is asymptotically normally distributed; furthermore
\begin{equation}\label{normal}
    \Vech \widehat{\Cov}(\overline{\mathbf{y}}_{_{ST}}) \build{\rightarrow}{d}{} \mathcal{N}_{k}
  \left(\E\left(\Vech \widehat{\Cov}(\overline{\mathbf{y}}_{_{ST}})\right),
  \Cov\left(\Vech \widehat{\Cov}(\overline{\mathbf{y}}_{_{ST}})\right)\right),
\end{equation}
where
\begin{equation}\label{ecyst}
    \E\left(\Vech \widehat{\Cov}(\overline{\mathbf{y}}_{_{ST}})\right) =  \sum_{h=1}^{H}\left(
     \frac{{{W_{h}}^{2}}}{n_{h}} - \frac{{W_{h}}}{N} \right) \frac{n_{h}}{n_{h}-1}\Vech\mathbf{S}_{h},
\end{equation}
\begin{eqnarray}
    \Cov\left(\Vech \widehat{\Cov}(\overline{\mathbf{y}}_{_{ST}})\right) \hspace{8cm}\nonumber \\
    \label{ccyst}
    \phantom{xx}=\sum_{h=1}^{H}\left(
    \frac{{{W_{h}}^{2}}}{n_{h}} - \frac{{W_{h}}}{N} \right)^{2} \frac{n_{h}}{(n_{h}-1)^{2}}
    \left(\mathbf{M}_{h}^{4} - \Vech \mathbf{S}_{h}\Vech' \mathbf{S}_{h}\right),
\end{eqnarray}
and
$$
  \mathbf{M}_{h}^{4} = \frac{1}{N_{h}}\mathbf{D}_{G}^{+}\left[\sum_{i = 1}^{N_{h}}
    (\mathbf{y}_{h i} - \overline{\mathbf{Y}}_{h})(\mathbf{y}_{h i} - \overline{\mathbf{Y}}_{h})'
    \otimes (\mathbf{y}_{h i} - \overline{\mathbf{Y}}_{h})(\mathbf{y}_{h i} - \overline{\mathbf{Y}}_{h})'
    \right]\mathbf{D}_{G}^{+'}.
$$
\end{thm}

Observe that the asymptotic means and covariance matrices of the asymptotically normality
distributions of $\overline{\mathbf{y}}_{h}$, $\Vech \mathbf{S}_{h}$, $\Vec
\widehat{\Cov}(\overline{\mathbf{y}}_{_{ST}})$ and $ \Vech
\widehat{\Cov}(\overline{\mathbf{y}}_{_{ST}})$ are in terms of the populations parameters
$\overline{\mathbf{Y}}_{h}$, $\Vech \mathbf{S}_{h}$, $\mathfrak{M}_{h}^{4}$ and
$\mathbf{M}_{h}^{4}$; then, from \citet[iv), pp. 388-389]{r:73}, approximations of asymptotic
distributions can be obtained using consistent estimators instead of population parametrers.
In what follows, the following substitutions are used:
\begin{equation}\label{sus}
    \overline{\mathbf{Y}}_{h} \rightarrow \overline{\mathbf{y}}_{h}, \qquad \Vech \mathbf{S}_{h}
  \rightarrow\Vech \mathbf{s}_{h},  \quad \mathfrak{M}_{h}^{4} \rightarrow \boldsymbol{\mathfrak{m}}_{h}^{4}
  \quad \mbox{ and } \quad \mathbf{M}_{h}^{4} \rightarrow \mathbf{m}_{h}^{4}
\end{equation}
where
$$
  \mathbf{m}_{h}^{4} = \frac{1}{n_{h}}\mathbf{D}_{G}^{+}\left[\sum_{i = 1}^{n_{h}}
    (\mathbf{y}_{h i} - \overline{\mathbf{y}}_{h})(\mathbf{y}_{h i} - \overline{\mathbf{y}}_{h})'
    \otimes (\mathbf{y}_{h i} - \overline{\mathbf{y}}_{h})(\mathbf{y}_{h i} - \overline{\mathbf{y}}_{h})'
    \right]\mathbf{D}_{G}^{+'},
$$
and
$$
  \boldsymbol{\mathfrak{m}}_{h}^{4} = \frac{1}{n_{h}}\left[\sum_{i = 1}^{n_{h}}
    (\mathbf{y}_{h i} - \overline{\mathbf{y}}_{h})(\mathbf{y}_{h i} - \overline{\mathbf{y}}_{h})'
    \otimes (\mathbf{y}_{h i} - \overline{\mathbf{y}}_{h})(\mathbf{y}_{h i} - \overline{\mathbf{y}}_{h})'
    \right].
$$

\section{\normalsize OPTIMUM ALLOCATION IN MULTIVARIATE STRATIFIED RANDOM SAMPLING VIA STOCHASTIC MATRIX MATHEMATICAL PROGRAMMING}\label{sec4}

When the variances are the objective functions, subject to certain cost function,  the optimum
allocation in multivariate stratified random sampling can be expressed as the following matrix
mathematical programming using a deterministic approach
\begin{equation}\label{om1}
  \begin{array}{c}
    \build{\min}{}{\mathbf{n}}{\widehat{\Cov}}(\overline{\mathbf{y}}_{_{ST}})\\
    \mbox{subject to}\\
    \mathbf{c}'\mathbf{n} + c_{0} = C  \\
    2\leq n_{h}\leq N_{h}, \ \ h=1,2,\dots, H\\
    n_{h}\in \mathbb{N},
  \end{array}
  \end{equation}
where $\mathbb{N}$ denotes the set of natural numbers. (\ref{om1}) has been studied in detail
by \citet{dgu:08}.

Observing that $\widehat{\Cov}(\overline{\mathbf{y}}_{_{ST}})$ is in terms of $s_{h_{jk}}$,
which are random variables, the optimum allocation of (\ref{om1}) via stochastic mathematical
programming can be stated as the following stochastic matrix mathematical programming, see
\citet{p:95} and \citet{sm:84},
\begin{equation}\label{smamp}
  \begin{array}{c}
    \build{\min}{}{\mathbf{n}}{\widehat{\Cov}}(\overline{\mathbf{y}}_{_{ST}})\\
    \mbox{subject to}\\
    \mathbf{c}'\mathbf{n} + c_{0} = C  \\
    2\leq n_{h}\leq N_{h}, \ \ h=1,2,\dots, H\\
    \Vech \widehat{\Cov}(\overline{\mathbf{y}}_{_{ST}}) \build{\rightarrow}{d}{} \mathcal{N}_{k}
  \left(\E\left(\Vech \widehat{\Cov}(\overline{\mathbf{y}}_{_{ST}})\right),
  \Cov\left(\Vech \widehat{\Cov}(\overline{\mathbf{y}}_{_{ST}})\right)\right)\\
    n_{h}\in \mathbb{N},
  \end{array}
\end{equation}
where $\E\left(\Vech \widehat{\Cov}(\overline{\mathbf{y}}_{_{ST}})\right)$ and $\Cov\left(\Vech
\widehat{\Cov}(\overline{\mathbf{y}}_{_{ST}})\right)$ are given by (\ref{ecyst}) and
(\ref{ccyst}) respectively.

Observe that $\widehat{\Cov}(\overline{\mathbf{y}}_{_{ST}})$ is an explicit function of
$\mathbf{n}$, and so it must be denoted as $\widehat{\Cov}(\overline{\mathbf{y}}_{_{ST}})
\equiv \widehat{\Cov}(\overline{\mathbf{y}}_{_{ST}}(\mathbf{n}))$. Also, assume that
$\widehat{\Cov}(\overline{\mathbf{y}}_{_{ST}}(\mathbf{n}))$ is a positive definite matrix for
all $\mathbf{n}$, $\widehat{\Cov}(\overline{\mathbf{y}}_{_{ST}}(\mathbf{n})) > \mathbf{0}$.
Now, let $\mathbf{n_{1}}$ and $\mathbf{n_{2}}$ be two possible values of the vector
$\mathbf{n}$ and, recall that, for $\mathbf{A}$ and $\mathbf{B}$ positive definite matrices,
$\mathbf{A} > \mathbf{B} \Leftrightarrow \mathbf{A} - \mathbf{B} > \mathbf{0}$.

Then, proceeding as \citet{dgu:08} the stochastic solution of (\ref{smamp}) is reduced to the
following stochastic uniobjective mathematical programming problem
\begin{equation}\label{smauo}
    \begin{array}{c}
    \build{\min}{}{\mathbf{n}}f\left(\widehat{\Cov}(\overline{\mathbf{y}}_{_{ST}})\right)\\
    \mbox{subject to}\\
    \mathbf{c}'\mathbf{n}+ c_{0}=C \\
    2\leq n_{h}\leq N_{h}, \ \ h=1,2,\dots, H\\
    \Vech \widehat{\Cov}(\overline{\mathbf{y}}_{_{ST}}) \build{\rightarrow}{d}{} \mathcal{N}_{k}
    \left(\E\left(\Vech \widehat{\Cov}(\overline{\mathbf{y}}_{_{ST}})\right),
    \Cov\left(\Vech \widehat{\Cov}(\overline{\mathbf{y}}_{_{ST}})\right)\right)\\
    n_{h}\in \mathbb{N},
  \end{array}
\end{equation}
where the function $f$ is such that: $f: \mathcal{S}\rightarrow \Re$,
\begin{equation}\label{citerio}
    \widehat{\Cov}(\overline{\mathbf{y}}_{_{ST}}(\mathbf{n_{1}})) <
    \widehat{\Cov}(\overline{\mathbf{y}}_{_{ST}}(\mathbf{n_{2}}))
    \Leftrightarrow
    f\left(\widehat{\Cov}(\overline{\mathbf{y}}_{_{ST}}(\mathbf{n_{1}}))\right)
    < f\left(\widehat{\Cov}(\overline{\mathbf{y}}_{_{ST}}(\mathbf{n_{2}}))\right).
\end{equation}
with $\widehat{\Cov}(\overline{\mathbf{y}}_{_{ST}}(\mathbf{n})) \in \mathcal{S}\subset
\Re^{G(G+1)/2}$ and $\mathcal{S}$ is the set of positive definite matrices.

Unfortunately or fortunately the function $f(\cdot)$ is not unique. Same alternatives for
$f\left(\widehat{\Cov}(\overline{\mathbf{y}}_{_{ST}}(\mathbf{n}))\right)$ are $ \tr
\left(\cdot\right)$,  $ \left|\cdot\right|$, $\lambda_{\max}\left(\cdot\right)$, where
$\lambda_{\max}$ is the maximum eigenvalue,  $ \lambda_{\min}\left(\cdot\right)$, where
$\lambda_{\min}$ is the minimum eigenvalue, $ \lambda_{j}\left(\cdot\right)$, where
$\lambda_{j}$ is the $j$-th eigenvalue, among others.

Note that (\ref{smauo}) is a stochastic uniobjective mathematical programming then, any
technique of stochastic uniobjective mathematical programming can be applied, for example:

Point $\mathbf{n} \in \mathbb{N}^{H}$ is the expected modified value solution to (\ref{smauo})
if it is an efficient solution in the \textbf{Pareto}\footnote{For the sampling context, observe
that in matrix mathematical programming problems, there rarely exists a point $\mathbf{n^{*}}$
which is considered as a minimum. Alternatively, it say that $f^{*}(\mathbf{x})$ is a
\textit{Pareto point} of $f(\mathbf{n}) = (f_{1}(\mathbf{n}), \dots, f_{G}(\mathbf{n}))'$, if
there is not other point $f^{1}(\mathbf{n})$ such that $f^{1}(\mathbf{n}) \leq
f^{*}(\mathbf{n})$, i.e. for all $j$, $f^{1}_{j}(\mathbf{n}) \leq f^{*}_{j}(\mathbf{n})$ and
$f^{1}(\mathbf{n}) \neq f^{*}(\mathbf{n})$.} sense to following deterministic uniobjetive
mathematical programming problem
\begin{equation}\label{emsma}
    \begin{array}{c}
    \build{\min}{}{\mathbf{n}} k_{1}\E\left(f\left(\widehat{\Cov}(\overline{\mathbf{y}}_{_{ST}})\right)\right)
            + k_{2}\sqrt{\Var\left(f\left(\widehat{\Cov}(\overline{\mathbf{y}}_{_{ST}})\right)\right)}\\
    \mbox{subject to}\\
    \mathbf{c}'\mathbf{n}+ c_{0}=C \\
    2\leq n_{h}\leq N_{h}, \ \ h=1,2,\dots, H\\
    n_{h}\in \mathbb{N},
  \end{array}
\end{equation}
Here $k_{1}$ and $k_{2}$ are non negative constants, and their values show the relative
importance of the expectation and the covariance matrix
$\widehat{\Cov}(\bar{\mathbf{y}}_{_{ST}})$. Some authors suggest that $k_{1} + k_{2} =1$, see
\citet[p. 599]{rao78}. Observe that if $k_{1}$ and $k_{2}$ are such that $k_{1} = 1$ and $k_{2}
= 0$ in (\ref{emsma}), the resulting method is known as the E-model. Alternatively, if $k_{1} =
0$ and $k_{2} = 1$, the method is called the V-model, see \citet{chc:63}, \citet{p:95} and
\citet{up:01}.

Alternatively, the point $\mathbf{n} \in \mathbb{N}^{H}$ is a minimum risk solution of the
aspiration level $\tau$ to the problem (\ref{smauo}) (also termed P-model, see \citet{chc:63})
if its is an efficient solution in the Pareto sense of the uniobjetive stochastic optimization
problem
\begin{equation}\label{mrsma}
    \begin{array}{c}
    \build{\min}{}{\mathbf{n}} \P\left(f\left(\widehat{\Cov}(\overline{\mathbf{y}}_{_{ST}})\right) \leq \tau\right)\\
    \mbox{subject to}\\
    \mathbf{c}'\mathbf{n}+ c_{0}=C \\
    2\leq n_{h}\leq N_{h}, \ \ h=1,2,\dots, H\\
    n_{h}\in \mathbb{N}.
  \end{array}
\end{equation}

In Section \ref{sec5} the solution is studied for the case when $f =
\tr\left(\widehat{\Cov}(\overline{\mathbf{y}}_{_{ST}})\right)$ and the case when $f =
\left|\widehat{\Cov}(\overline{\mathbf{y}}_{_{ST}})\right|$. These solutions are implemented in the context of
problems (\ref{emsma}) and (\ref{mrsma}).

Finally, note that so far, the cost constraint $\displaystyle\sum_{h=1}^{H}c_{h}n_{h}+
c_{0}=C$  has been used in every stochastic mathematical programming method. However, in diverse situations,
this cost restriction could represent existing restrictions on the availability of man-hours for
carrying out a survey, or restrictions on the total available time for performing the survey,
etc. These limitations can be established by using the following constraint, see \cite{ad81}:
$$
  \sum_{h=1}^{H}n_{h} = n.
$$

\section{\normalsize APPLICATION}\label{sec5}

The input information was taken from \citet{aa81} in which they describe a forest survey conducted in Humbolt
County, California. The population was subdivided into nine strata on the basis of the timber
volume per unit area, as determined from aerial photographs. The two variables included in this
example are the basal area (BA)\footnote{In forestry terminology, `Basal area' is the area of a
plant perpendicular to the longitudinal axis of a tree at 4.5 feet above ground.} in square
feet, and the net volume in cubic feet (Vol.), both expressed on a per acre basis. The
variances, covariances and the number of units within stratum $h$ are listed in Table 1.

\begin{table}
\caption{\small Variances, covariances and the number of units within each stratum}
\begin{center}
\begin{footnotesize}
\begin{tabular}{ c r r r r }
\hline\hline
\multicolumn{2}{c}{} & \multicolumn{2}{c}{Variance} \\
\cline{3-4}
Stratum & $N_{h}$ & \hspace{.5cm} BA \hspace{.5cm} & \hspace{.5cm} Vol. \hspace{.5cm} & \hspace{.5cm}Covariance \\
\hline\hline
1 & 11 131 & 1 557 & 554 830 & 28 980 \\
2 & 65 857 & 3 575 & 1 430 600 & 61 591\\
3 & 106 936 & 3 163 & 1 997 100 & 72 369 \\
4 & 72 872 & 6 095 & 5 587 900 & 166 120\\
5 & 78 260 & 10 470 & 10 603 000 & 293 960 \\
6 & 51 401 & 8 406 & 15 828 000 & 357 300\\
7 & 24 050 & 20 115 & 26 643 000 & 663 300 \\
8 & 46 113 & 9 718 & 13 603 000 & 346 810\\
9 & 102 985 & 2 478 & 1 061 800 & 39 872 \\
\hline\hline
\end{tabular}
\end{footnotesize}
\end{center}
\end{table}

For this example, the matrix optimisation problem under approach (\ref{smauo}) is
\begin{equation}\label{ej}
  \begin{array}{c}
  \build{\min}{}{\mathbf{n}}
    f\left(%
    \begin{array}{c c}
      \widehat{\Var}(\overline{y}_{_{ST}}^{1}) & \widehat{\Cov}(\overline{y}_{_{ST}}^{1}, \overline{y}_{_{ST}}^{2})\\
      \widehat{\Cov}(\overline{y}_{_{ST}}^{2}, \overline{y}_{_{ST}}^{1}) & \widehat{\Var}(\overline{y}_{_{ST}}^{2}) \\
    \end{array}%
    \right)\\
    \mbox{subject to}\\
    \displaystyle\sum_{h=1}^{9}n_{h}=1000 \\
    2\leq n_{h}\leq N_{h}, \ \ h=1,\dots, 9\\
    \Vech \widehat{\Cov}(\overline{\mathbf{y}}_{_{ST}}) \build{\rightarrow}{d}{} \mathcal{N}_{3}
  \left(\E\left(\Vech \widehat{\Cov}(\overline{\mathbf{y}}_{_{ST}})\right),
  \Cov\left(\Vech \widehat{\Cov}(\overline{\mathbf{y}}_{_{ST}})\right)\right)\\
    n_{h}\in \mathbb{N}.
  \end{array}
\end{equation}

\subsection{\normalsize Solution when $f(\cdot) \equiv \tr(\cdot)$}

Note that by (\ref{normal}), (\ref{ecyst}) and (\ref{ccyst})
$$
  \tr \Cov\left(\overline{\mathbf{y}}_{ST}\right) \sim \mathcal{N}\left(\E \left (\tr \Cov\left(\overline{\mathbf{y}}_{ST}\right)\right),
  \Var\left(\tr \Cov\left(\overline{\mathbf{y}}_{ST}\right)\right)\right)
$$
where
$$
    \E\left(\tr \widehat{\Cov}(\overline{\mathbf{y}}_{_{ST}})\right) =  \sum_{j=1}^{G}\sum_{h=1}^{H}\left(
     \frac{{{W_{h}}^{2}}}{n_{h}} - \frac{{W_{h}}}{N} \right) \frac{n_{h}}{n_{h}-1}S_{h_{j}}^{2},
$$
$$
    \Var\left(\tr \widehat{\Cov}(\overline{\mathbf{y}}_{_{ST}})\right)
    =\sum_{j=1}^{G}\sum_{h=1}^{H}\left(
    \frac{{{W_{h}}^{2}}}{n_{h}} - \frac{{W_{h}}}{N} \right)^{2} \frac{n_{h}}{(n_{h}-1)^{2}}
    \left(m_{h_{j}}^{4} - (S_{h_{j}}^{2})^{2}\right),
$$
and
$$
  m_{h_{j}}^{4} = \frac{1}{N_{h}}\left[\sum_{i = 1}^{N_{h}} \left(y_{h i}^{j} - \overline{Y}_{h}^{j}\right)^{4} \right].
$$
Therefore, considering the substitutions (\ref{sus}), the equivalent deterministic uniobjetive
mathematical programming problem to stochastic mathematical programming (\ref{ej}) via the
modified $E$-model is
$$
    \begin{array}{c}
    \build{\min}{}{\mathbf{n}} k_{1}\widehat{\E}\left(\tr\widehat{\Cov}(\overline{\mathbf{y}}_{_{ST}})\right)
            + k_{2}\sqrt{\widehat{\Var}\left(\tr\widehat{\Cov}(\overline{\mathbf{y}}_{_{ST}})\right)}\\
    \mbox{subject to}\\
    \displaystyle\sum_{h=1}^{9}n_{h}=1000 \\
    2\leq n_{h}\leq N_{h}, \ \ h=1,2,\dots, 9\\
    n_{h}\in \mathbb{N},
  \end{array}
$$
where
\begin{equation}\label{esp}
    \widehat{\E}\left(\tr \widehat{\Cov}(\overline{\mathbf{y}}_{_{ST}})\right) =  \sum_{j=1}^{2}\sum_{h=1}^{9}\left(
     \frac{{{W_{h}}^{2}}}{n_{h}} - \frac{{W_{h}}}{N} \right) \frac{n_{h}}{n_{h}-1}s_{h_{j}}^{2},
\end{equation}
\begin{equation}\label{varia}
    \widehat{\Var}\left(\tr \widehat{\Cov}(\overline{\mathbf{y}}_{_{ST}})\right)
    =\sum_{j=1}^{2}\sum_{h=1}^{9}\left(
    \frac{{{W_{h}}^{2}}}{n_{h}} - \frac{{W_{h}}}{N} \right)^{2} \frac{n_{h}}{(n_{h}-1)^{2}}
    \left(\mathfrak{m}_{h_{j}}^{4} - (s_{h_{j}}^{2})^{2}\right),
\end{equation}
and
\begin{equation}\label{impo}
    \mathfrak{m}_{h_{j}}^{4} = \frac{1}{n_{h}}\left[\sum_{i = 1}^{n_{h}} \left(y_{h i}^{j} - \overline{y}_{h}^{j}\right)^{4} \right].
\end{equation}
\begin{rem}\label{impo1}
Observe that the estimators $\overline{y}_{h}^{j}$, $s_{h_{j}}^{2}$ and $m_{h_{j}}^{4}$ of
$\overline{Y}_{h}^{j}$, $S_{h_{j}}^{2}$ and $M_{h_{j}}^{4}$  are initially obtained as
\begin{enumerate}[i)]
  \item a consequence of a pilot study (or preliminary sample) or
  \item using the corresponding values of the estimators of another variable $X$ correlated to the variable $Y$.
\end{enumerate}
It is important to have this in mind in the the minimisation step, because for example,
the $n_{h}$'s that appear in expression (\ref{impo}), are the fixed $n_{h}$'s values used in the
pilot study. Same comment for the expressions of the estimator $\overline{y}_{h}^{j}$ and
$s_{h_{j}}^{2}$. While the $n_{h}$'s that appear in expressions (\ref{esp}) and (\ref{varia})
are the decision variables. \qed
\end{rem}

Similarly, proceeding as in \citet{dgrc:00}, and noting that, if $\Phi$ denotes the
distribution function of the standard Normal distribution, the objective function in (\ref{ej})
with $f(\cdot) \equiv \tr(\cdot)$ can be written as
$$
  \build{\min}{}{\mathbf{n}}\quad \Phi\left( \frac{\tau
  -\widehat{\E}\left(\tr\widehat{\Cov}(\overline{\mathbf{y}}_{_{ST}})\right)}
  {\sqrt{\widehat{\Var}\left(\tr\widehat{\Cov}(\overline{\mathbf{y}}_{_{ST}})\right)}} \right ).
$$
In this way, since minimising the  monotonically increasing distribution function is equivalent
to minimising the value of the associated random variable, the equivalent deterministic problem
to the stochastic mathematical programming (\ref{ej}) via the $P$-model is
$$
    \begin{array}{c}
    \build{\min}{}{\mathbf{n}} \displaystyle\frac{\tau -\widehat{\E}\left(\tr\widehat{\Cov}(\overline{\mathbf{y}}_{_{ST}})\right)}
     {\sqrt{\widehat{\Var}\left(\tr\widehat{\Cov}(\overline{\mathbf{y}}_{_{ST}})\right)}}\\
    \mbox{subject to}\\
    \displaystyle\sum_{h=1}^{9}n_{h}=1000 \\
    2\leq n_{h}\leq N_{h}, \ \ h=1,2,\dots, 9\\
    n_{h}\in \mathbb{N},
  \end{array}
$$

\begin{rem}

When $f(\cdot) \equiv |\cdot|$, this approach consider the following alternative stochastic
matrix mathematical programming problem
\begin{equation}\label{det1}
  \begin{array}{c}
    \build{\min}{}{\mathbf{n}}\widehat{\Cov}(\overline{\mathbf{y}}_{_{ST}})\\
    \mbox{subject to}\\
    \displaystyle\sum_{h=1}^{9}n_{h}=1000 \\
    2\leq n_{h}\leq N_{h}, \ \ h=1,2,\dots, 9\\
    \Vech \widehat{\Cov}(\overline{\mathbf{y}}_{_{ST}}) \build{\rightarrow}{d}{} \mathcal{N}_{2 \times 2}
    \left(\Vech \mathbf{0}_{2 \times 2},
    \Cov\left(\Vech\widehat{\Cov}(\overline{\mathbf{y}}_{_{ST}})\right)\right)\\
    n_{h}\in \mathbb{N},
  \end{array}
\end{equation}
where $\widehat{\Cov}(\overline{\mathbf{y}}_{_{ST}})
   = \Vech^{-1}\left[\Vech\widehat{\Cov}(\overline{\mathbf{y}}_{_{ST}})
  - \E\left(\Vech\widehat{\Cov}(\overline{\mathbf{y}}_{_{ST}})\right)\right]
$
and $\Vech^{-1}$ is the inverse function of function $\Vech$.

In this way (\ref{mrsma}) is
\begin{equation}\label{det2}
  \begin{array}{c}
    \build{\min}{}{\mathbf{n}}\left|\widehat{\Cov}(\overline{\mathbf{y}}_{_{ST}})\right|\\
    \mbox{subject to}\\
    \displaystyle\sum_{h=1}^{9}n_{h}=1000 \\
    2\leq n_{h}\leq N_{h}, \ \ h=1,2,\dots, 9\\
    \Vech \widehat{\Cov}(\overline{\mathbf{y}}_{_{ST}}) \build{\rightarrow}{d}{} \mathcal{N}_{2 \times 2}
    \left(\Vech \mathbf{0}_{2 \times 2},
    \Cov\left(\Vech\widehat{\Cov}(\overline{\mathbf{y}}_{_{ST}})\right)\right)\\
    n_{h}\in \mathbb{N},
  \end{array}
\end{equation}
Thus, taking into account the substitutions (\ref{sus}), the equivalent deterministic uniobjetive
mathematical programming problem to the stochastic mathematical programming (\ref{det2}) via the
modified $E$-model is
$$
    \begin{array}{c}
    \build{\min}{}{\mathbf{n}}
    k_{1}\widehat{\E}\left(\left|\widehat{\Cov}(\overline{\mathbf{y}}_{_{ST}})\right|\right)
    + k_{2}\sqrt{\widehat{\Var}\left(\left|\widehat{\Cov}(\overline{\mathbf{y}}_{_{ST}})\right|\right)}\\
    \mbox{subject to}\\
    \displaystyle\sum_{h=1}^{9}n_{h}=1000 \\
    2\leq n_{h}\leq N_{h}, \ \ h=1,2,\dots, 9\\
    n_{h}\in \mathbb{N},
  \end{array}
$$
where for $G = 2$ and assuming that
$\widehat{\Cov}\left(\Vech\widehat{\Cov}(\overline{\mathbf{y}}_{_{ST}})\right)$ is such that
$$
  \widehat{\Cov}\left(\Vech\widehat{\Cov}(\overline{\mathbf{y}}_{_{ST}})\right) = \mathbf{B}\otimes
  \mathbf{B},
$$
it is obtained that, see \citet{dc:00},
$$
    \widehat{\E}\left(\left| \widehat{\Cov}(\overline{\mathbf{y}}_{_{ST}})\right|\right) =
    |\mathbf{N}|^{1/4} \frac{(-1)}{\sqrt{\pi}}\left(\Gamma[1/2]-\Gamma[3/2]\right),
$$
and $\widehat{\Var}\left(\left|\widehat{\Cov}(\overline{\mathbf{y}}_{_{ST}})\right|\right)$
is%
$$
        =|\mathbf{N}|^{1/2} \left[\frac{2}{\sqrt{\pi}}\left(\Gamma[1/2]-\Gamma[3/2] + \frac{\Gamma[5/2]}{2}\right)
    -\frac{1}{\pi}\left(\Gamma[1/2]-\Gamma[3/2]\right)^{2}\right],
$$
where $\Gamma[\cdot]$ denotes the gamma function,
$$
  \mathbf{N} = \sum_{h=1}^{H}\left(
    \frac{{{W_{h}}^{2}}}{n_{h}} - \frac{{W_{h}}}{N} \right)^{2} \frac{n_{h}}{(n_{h}-1)^{2}}
    \left(\mathbf{\mathfrak{m}}_{h}^{4} - \Vec \mathbf{s}_{h}\Vec' \mathbf{s}_{h}\right)
$$
and
$$
  \boldsymbol{\mathfrak{m}}_{h}^{4} = \frac{1}{n_{h}}\left[\sum_{i = 1}^{n_{h}}
    (\mathbf{y}_{h i} - \overline{\mathbf{y}}_{h})(\mathbf{y}_{h i} - \overline{\mathbf{y}}_{h})'
    \otimes (\mathbf{y}_{h i} - \overline{\mathbf{y}}_{h})(\mathbf{y}_{h i} - \overline{\mathbf{y}}_{h})'
    \right],
$$
see Remark \ref{impo1}.

Similarly, considering (\ref{det1}) and that $f(\cdot) \equiv |\cdot|$, (\ref{mrsma}) is
restated as
$$
    \begin{array}{c}
    \build{\min}{}{\mathbf{n}} \P\left(\left|\widehat{\Cov}(\overline{\mathbf{y}}_{_{ST}})\right| \leq \tau\right)\\
    \mbox{subject to}\\
    \displaystyle\sum_{h=1}^{9}n_{h}=1000 \\
    2\leq n_{h}\leq N_{h}, \ \ h=1,2,\dots, 9\\
    \Vech \widehat{\Cov}(\overline{\mathbf{y}}_{_{ST}}) \build{\rightarrow}{d}{} \mathcal{N}_{2 \times 2}
    \left(\Vech \mathbf{0}_{2 \times 2},
    \Cov\left(\Vech\widehat{\Cov}(\overline{\mathbf{y}}_{_{ST}})\right)\right)\\
    n_{h}\in \mathbb{N}
  \end{array}
$$
Then, if $\Psi$ denotes the distribution function of the determinant of
$\widehat{\Cov}(\overline{\mathbf{y}}_{_{ST}})$,  the equivalent deterministic problem to
the stochastic mathematical programming (\ref{ej}) via the $P$-model is
$$
    \begin{array}{c}
    \build{\min}{}{\mathbf{n}}\quad \tau |\mathbf{N}|^{1/4}\\
    \mbox{subject to}\\
    \displaystyle\sum_{h=1}^{9}n_{h}=1000 \\
    2\leq n_{h}\leq N_{h}, \ \ h=1,2,\dots, 9\\
    n_{h}\in \mathbb{N},
  \end{array}
$$
where the density of $Z = \widehat{\Cov}(\overline{\mathbf{y}}_{_{ST}})$ is, see \citet{dc:00}
$$
  \frac{dG(z)}{dz} = g_{_{Z}}(z) = \frac{1}{\sqrt{2}} \exp (z)\left[1 - \erf\left(\sqrt{2z}\right)\right], \quad
  z\geq 0,
$$
where $\erf(\cdot)$ is the usual error function defined as
$$
  \erf(x) = \frac{2}{\sqrt{\pi}}\int_{0}^{x} \exp(-t^{2})dt.
$$ \qed
\end{rem}

Table 2 shows the optimisation solutions obtained by some of the methods described in Section
\ref{sec4}. Specifically, the solution is presented for the case when the value function is
defined as the trace function, $f(\cdot) = \tr(\cdot)$ and for the following stochastic solutions:
Modified $E-$model, $E-$model, $V-$model and the $P-$model. Also, the
optimum allocation is included for each characteristic, BA and Vol (the first two rows in Table 2). The
last two columns show the minimum values of the individual variances for the respective optimum
allocations identified by each method. The results were computed using the commercial software
Hyper LINGO/PC, release 6.0, see \citet{w95}. The default optimisation methods used by LINGO to
solve the nonlinear integer optimisation programs are Generalised Reduced Gradient (GRG) and
branch-and-bound methods, see \citet{bss06}. Some technical details of the computations are the
following: the maximum number of iterations of the methods presented in Table 2 was 2279
(modified $E$-model) and the mean execution time for all the programs was 4 seconds.  Finally,
note that the greatest discrepancy found by the different methods among the sizes of the strata
occurred under $P$-model. Beyond doubt, this is a consequence of the election of the corresponding
value of $\tau$ needed for the $P$-model approach.

\begin{table}
\caption{\small Sample sizes and estimator of variances for the different allocations
calculated}
\begin{center}
\begin{minipage}[t]{400pt}
\begin{scriptsize}
\begin{tabular}{ c  c  c  c  c  c  c c c c c c}
\hline\hline Allocation\footnote{The estimated fourth moment $m_{h_{j}}^{4}$ were simulated.} &
$n_{1}$ & $n_{2}$ & $n_{3}$ & $n_{4}$ & $n_{5}$ & $n_{6}$ & $n_{7}$ & $n_{8}$ & $n_{9}$ &
$\widehat{\Var}(\overline{y}_{_{ST}}^{1})$ &
$\widehat{\Var}(\overline{y}_{_{ST}}^{2})$ \\
\hline\hline
BA & 10 & 94 & 144 & 136 & 191 & 113 & 81 & 109 & 122 & 5.591 & 5441.105\\
Vol & 7 & 62 & 119 & 136 & 200 & 161 & 98 & 134 & 83 & 5.953 & 5139.531\\
$\boldsymbol{\tr \widehat{\Cov}(\overline{\mathbf{y}}_{_{ST}})}$
& &  &  &  &  & &&  &  & & \\
\begin{tabular}{c}
  Modified \\
  $E$-model \\
\end{tabular}
& 8 & 46 & 77 & 119 & 191 & 191 & 158 & 161 & 49 & 7.312 & 5593.494\\
$E$-model\footnote{Where $k_{1} = k_{2} = 0.5$.} & 7 & 63 & 119 & 135 & 200 & 160 & 98 & 134 & 84 & 5.937 & 5139.645\\
$V$-model & 8 & 46 & 77 & 119 & 191 & 191 & 158 & 161 & 49 & 7.312 & 5593.494\\
$P$-model\footnote{Where $\tau = 6000$.} & 632 & 9 & 117 & 29 & 46 & 54 & 52 & 49 & 7 & 29.746 & 20820.660\\
\hline\hline
\end{tabular}
\end{scriptsize}
\end{minipage}
\end{center}
\end{table}

\section*{\normalsize CONCLUSIONS}

It is difficult to suggest general rules for the selection of a method in stochastic matrix
mathematical programming (\ref{smamp}). These conclusions are sustained in several regards, for
example: potentiality, there is an infinite number of possible definitions of the value
function $f(\cdot)$; furthermore, the value function approach is not the unique way to restate
(\ref{smamp}); exist many ways to solve (\ref{smamp}) from a stochastic point of view.
We believe that this responsibility lies with the person skilled in the particular field and in
his/her capacity of discern which function or approach that better reflects and meets the
objectives of the study.

In this paper, the problem of optimal allocation in multivariate stratified sampling
was considered. In all sample size problems there is always uncertainty regarding the population parameters
and in this work, this uncertainty was incorporated via a stochastic matrix mathematical solution.

\section*{\normalsize ACKNOWLEDGMENTS}

This research work was partially supported by IDI-Spain, Grants No. FQM2006-2271 and
MTM2008-05785, supported also by CONACYT Grant CB2008 Ref. 105657. This paper was written during J. A. D\'{\i}az-Garc\'{\i}a's stay as a visiting
professor at the Department of Probability Statistics of the Center of Mathematical Research,
Guanajuato, M\'exico.

\end{document}